\documentclass[a4paper,reqno,12pt]{amsart}
\usepackage[margin=1.1in]{geometry}
\usepackage{indentfirst}
\usepackage{mathrsfs}
\usepackage{syntonly}
\usepackage{amsmath}
\usepackage{amsthm}
\usepackage{amsfonts}
\usepackage{amssymb}
\usepackage{latexsym}
\usepackage{amscd,amssymb,amsopn,amsmath,amsthm,graphics,amsfonts,mathrsfs,accents,enumerate,verbatim,calc}
\usepackage[dvips]{graphicx}
\usepackage[colorlinks=true,linkcolor=blue,citecolor=blue]{hyperref}
\input xy
\xyoption{all}
\def\del{\delta}
\usepackage{color}

\numberwithin{equation}{section}
\theoremstyle{plain}

\newtheorem{thm}{Theorem}[section]
\newtheorem{cor}[thm]{Corollary}
\newtheorem{lem}[thm]{Lemma}

\newtheorem{defn}[thm]{Definition}
\newtheorem{exm}[thm]{Example}
\theoremstyle{remark}
\newtheorem{rem}[thm]{Remark}

\newcommand{\coresdim}{{\rm coresdim}}
\newcommand{\resdim}{{\rm resdim}}
\def\C{\mathscr{C}}
\def\X{\mathscr{X}}
\def\Y{\mathscr{Y}}

\def\I{\mathcal{I}}
\def\P{\mathcal{P}}

\def\Hom{\mbox{Hom}}
\def\E{\mathbb{E}}
\def\s{\mathfrak{s}}
\usepackage{bm}
\usepackage{tikz}
\usepackage{enumitem}
\setenumerate[1]{itemsep=4pt}
\begin{document}
\title[$n$-cotorsion pairs and $(n+1)$-cluster tilting subcategories]{On the relation between $n$-cotorsion pairs and \\[2mm]$(n+1)$-cluster tilting subcategories}\footnote{This work is supported by the Hunan Provincial Natural Science Foundation of China (Grants No. 2018JJ3205) and the NSF of China (Grants No. 11671221).}

\author[P. Zhou]{Panyue Zhou}
\address{College of Mathematics, Hunan Institute of Science and Technology, 414006 Yueyang, Hunan, People's Republic of China.}
\email{panyuezhou@163.com}

\baselineskip=16pt

\begin{abstract}
A notion of $n$-cotorsion pairs in an extriangulated category with enough projectives and enough injectives is defined in this article. We show that there exists a one-to-one correspondence between $n$-cotorsion pairs and $(n+1)$-cluster tilting subcategories. As an application, this result generalizes the work by Huerta, Mendoza and P\'{e}rez in an abelian case. Finally, we give some examples illustrating our main result.
\end{abstract}
\maketitle

\medskip

\section{Introduction}
Motivated by some properties satisfied by Gorenstein projective and Gorenstein injective modules over an Iwanaga-Gorenstein ring, Huerta, Mendoza and P\'{e}rez \cite[Definition 2.2]{HMP} introduced the notion of left and right $n$-cotorsion pairs in an abelian category $\C$. Two classes $\X$ and $\Y$ of objects of $\C$ form a \emph{left $n$-cotorsion pair} $(\X,\Y)$ in $\C$ if the orthogonality relation $\mathsf{Ext}^k_{\C}(\X,\Y) = 0$ is satisfied for any $1 \leq k \leq n$, and if every object of $\C$ has a resolution by objects in $\X$ whose syzygies have $\Y$-resolution dimension at most $n-1$.
Dually we can define the notion of a right $n$-cotorsion pair. If $(\X,\Y)$ is both a left and right $n$-cotorsion pair in $\C$, we call $(\X,\Y)$ an \emph{$n$-cotorsion pair}.
This concept generalises the notion of complete cotorsion pairs.  They also showed the following.
\begin{thm}\emph{\cite[Theorem 5.26]{HMP}}
Let $\C$ be an abelian category with enough projectives and enough injectives.
Then for any subcategory $\X$\hspace{-0.8mm} of\hspace{0.8mm} $\C$ and any integer $n\geq 1$, the following statements are equivalent:
\begin{enumerate}
\item[{\rm (1)}] $(\X,\X)$ is an $n$-cotorsion pair in $\C$.

\item[{\rm (2)}] $\X$ is an $(n+1)$-cluster tilting subcategory of $\C$.
\end{enumerate}
\end{thm}
Recently, the notion of extriangulated categories was introduced in \cite{NP} as a simultaneous generalization of exact categories (abelian categories are also exact categories) and triangulated categories. Exact categories and extension closed subcategories of a
triangulated category are extriangulated categories, while there are some other examples of extriangulated categories which are neither exact nor triangulated, see \cite[Proposition 3.30]{NP}, \cite[Example 4.14]{ZZ1} and \cite[Remark 3.3]{HZZ}.  Motivated by this idea, we introduce $n$-cotorsion pairs in an extriangulated category with enough projectives and enough injectives, for more details, see Definition \ref{d1}.
Our main result is the following.
\begin{thm}\emph{(See Theorem \ref{main})}
Let $\C$ be an extriangulated category with enough projectives and enough injectives.
Then for any subcategory $\X$\hspace{-1.2mm} of\hspace{1mm} $\C$ and any integer $n\geq 1$, the following statements are equivalent:
\begin{enumerate}
\item[{\rm (1)}] $(\X,\X)$ is an $n$-cotorsion pair in $\C$.

\item[{\rm (2)}] $\X$ is an $(n+1)$-cluster tilting subcategory of $\C$.
\end{enumerate}
\end{thm}
Since any abelan category is an extriangulated category, our main result generalizes
the work by Huerta, Mendoza and P\'{e}rez.  Note that any triangulated category can be viewed as an extriangulated category with enough projectives and enough injectives.
Our main result seems to be new phenomenon when it is applied to triangulated categories.

\begin{cor}
Let $\C$ be a triangulated category. Then for any subcategory $\X$\hspace{-1.2mm} of\hspace{1mm} $\C$  and any integer $n\geq 1$, the following statements are equivalent:
\begin{itemize}
\item[{\rm (1)}] $(\X,\X)$ is an $n$-cotorsion pair in $\C$.

\item[{\rm (2)}] $\X$ is an $(n+1)$-cluster tilting subcategory of $\C$.
\end{itemize}
\end{cor}

This article is organized as follows. In section 2, we recall some definitions and useful facts on extriangulated categories. In section 3, we prove our main result and give some examples illustrating it.

\section{Preliminaries}
Let us briefly recall some definitions and basic properties of extriangulated categories from \cite{NP}.
We omit some details here, but the reader can find
them in \cite{NP}.

Let $\C$ be an additive category equipped with an additive bifunctor
$$\mathbb{E}: \C^{\rm op}\times \C\rightarrow {\rm Ab},$$
where ${\rm Ab}$ is the category of abelian groups. For any objects $A, C\in\C$, an element $\delta\in \mathbb{E}(C,A)$ is called an $\mathbb{E}$-extension.
Let $\mathfrak{s}$ be a correspondence which associates an equivalence class $$\mathfrak{s}(\delta)=\xymatrix@C=0.8cm{[A\ar[r]^x
 &B\ar[r]^y&C]}$$ to any $\mathbb{E}$-extension $\delta\in\mathbb{E}(C, A)$. This $\mathfrak{s}$ is called a {\it realization} of $\mathbb{E}$, if it makes the diagrams in \cite[Definition 2.9]{NP} commutative.
 A triplet $(\C, \mathbb{E}, \mathfrak{s})$ is called an {\it extriangulated category} if it satisfies the following conditions.
\begin{enumerate}
\item $\mathbb{E}\colon\C^{\rm op}\times \C\rightarrow \rm{Ab}$ is an additive bifunctor.

\item $\mathfrak{s}$ is an additive realization of $\mathbb{E}$.

\item $\mathbb{E}$ and $\mathfrak{s}$  satisfy the compatibility conditions in \cite[Definition 2.12]{NP}.

 \end{enumerate}

We will use the following terminology.
\begin{defn}{\cite{NP}}
Let $\C$ be an extriangulated category.
\begin{enumerate}
\item[{\rm (1)}] A sequence $A\xrightarrow{~x~}B\xrightarrow{~y~}C$ is called a {\it conflation} if it realizes some $\E$-extension $\del\in\E(C,A)$. In this case, $x$ is called an {\it inflation} and $y$ is called a {\it deflation}.

\item[{\rm (2)}] If a conflation  $A\xrightarrow{~x~}B\xrightarrow{~y~}C$ realizes $\delta\in\mathbb{E}(C,A)$, we call the pair $( A\xrightarrow{~x~}B\xrightarrow{~y~}C,\delta)$ an {\it $\E$-triangle}, and write it in the following way.
$$A\overset{x}{\longrightarrow}B\overset{y}{\longrightarrow}C\overset{\delta}{\dashrightarrow}$$
We usually do not write this $``\delta"$ if it is not used in the argument.
\item[{\rm (3)}] Let $A\overset{x}{\longrightarrow}B\overset{y}{\longrightarrow}C\overset{\delta}{\dashrightarrow}$ and $A^{\prime}\overset{x^{\prime}}{\longrightarrow}B^{\prime}\overset{y^{\prime}}{\longrightarrow}C^{\prime}\overset{\delta^{\prime}}{\dashrightarrow}$ be any pair of $\E$-triangles. If a triplet $(a,b,c)$ realizes $(a,c)\colon\delta\to\delta^{\prime}$, then we write it as
$$\xymatrix{
A \ar[r]^x \ar[d]^a & B\ar[r]^y \ar[d]^{b} & C\ar@{-->}[r]^{\del}\ar[d]^c&\\
A'\ar[r]^{x'} & B' \ar[r]^{y'} & C'\ar@{-->}[r]^{\del'} &}$$
and call $(a,b,c)$ a {\it morphism of $\E$-triangles}.

\item[{\rm (4)}] An object $P\in\C$ is called {\it projective} if
for any $\E$-triangle $A\overset{x}{\longrightarrow}B\overset{y}{\longrightarrow}C\overset{\delta}{\dashrightarrow}$ and any morphism $c\in\C(P,C)$, there exists $b\in\C(P,B)$ satisfying $yb=c$.
We denote the subcategory of projective objects by $\mathcal P\subseteq\C$. Dually, the subcategory of injective objects is denoted by $\I\subseteq\C$.

\item[{\rm (5)}] We say that $\C$ {\it has enough projective objects} if
for any object $C\in\C$, there exists an $\E$-triangle
$A\overset{x}{\longrightarrow}P\overset{y}{\longrightarrow}C\overset{\delta}{\dashrightarrow}$
satisfying $P\in\mathcal P$. Dually we can define $\C$ {\it has enough injective objects}.
\end{enumerate}
\end{defn}

Let $\C$ be extriangulated category with enough projectives and enough injectives, and $\X$ a subcategory of $\C$.
We denote $\Omega\X={\rm CoCone}(\P,\X)$, that is to say, $\Omega\X$ is the subcategory of $\C$ consisting of
objects $\Omega X$ such that there exists an $\E$-triangle:
$$\Omega X\overset{a}{\longrightarrow}P\overset{b}{\longrightarrow}X\overset{}{\dashrightarrow},$$
with $P\in\P$ and $X\in\X$. We call $\Omega$ the \emph{syzygy} of $\X$.
Dually we define the \emph{cosyzygy} of $\X$ by $\Sigma\X={\rm Cone}(\X,\I)$.
Namely, $\Sigma\X$ is the subcategory of $\C$ consisting of
objects $\Sigma X$ such that there exists an $\E$-triangle:
$$X\overset{c}{\longrightarrow}I\overset{d}{\longrightarrow}\Sigma X\overset{}{\dashrightarrow},$$
with $I\in\I$ and $X\in\X$. For more details, see \cite[Definition 4.2 and Proposition 4.3]{LN}.

For a subcategory $\X\subseteq\C$, put $\Omega^0\X=\X$, and define $\Omega^k\X$ for $k>0$ inductively by
$$\Omega^k\X=\Omega(\Omega^{k-1}\X)={\rm CoCone}(\P,\Omega^{k-1}\X).$$
We call $\Omega^k\X$ the $k$-th syzygy of $\X$. Dually we define the $k$-th cosyzygy $\Sigma^k\X$ by
$\Sigma^0\X=\X$ and $\Sigma^k\X={\rm Cone}(\Sigma^{k-1}\X,\I)$ for $k>0$.

Liu and Nakaoka \cite{LN} defined higher extension groups in an extriangulated category with have
enough projectives and enough injectives as $\E(A,\Sigma^kB)\simeq\E(\Omega^kA,B)$ for $k\geq0$.
For convenience, we denote $\E(A,\Sigma^kB)\simeq\E(\Omega^kA,B)$ by $\E^{k+1}(A,B)$ for $k\geq0$.
They proved the following.
\begin{lem}
Let $\C$ be an extriangulated category with enough projectives and enough injectives. Assume that
$$\xymatrix{A\ar[r]^{f}& B\ar[r]^{g}&C\ar@{-->}[r]^{\delta}&}$$
is an $\E$-triangle in $\C$. Then for any object $X\in\C$ and $k\geq1$, we have the following exact sequences:
$$\cdots\xrightarrow{}\E^k(X,A)\xrightarrow{}\E^k(X,B)\xrightarrow{}\E^k(X,C)\xrightarrow{}\E^{k+1}(X,A)\xrightarrow{}\E^{k+1}(X,B)\xrightarrow{}\cdots;$$
$$\cdots\xrightarrow{}\E^k(C,X)\xrightarrow{}\E^k(B,X)\xrightarrow{}\E^k(A,X)\xrightarrow{}\E^{k+1}(C,X)\xrightarrow{}\E^{k+1}(B,X)\xrightarrow{}\cdots.$$
\end{lem}

As a higher version cluster tiling subcategories of extriangulated categories \cite[Definition 4.1]{CZZ}.
Liu and Nakaoka \cite[Definition 5.3]{LN} introduced the notion of $n$-cluster tiling subcategories of extriangulated categories. This definition generalizes Iyama's definition \cite[Definition 1.1]{I} in abelian case.

\begin{defn}\cite[Definition 5.3]{LN}
Let $\C$ be an extriangulated category with enough projectives and enough injectives.
A subcategory $\X\subseteq\C$ is called \emph{$n$-cluster tilting}, if it satisfies the
following conditions.
\begin{enumerate}
\item[{\rm (1)}] $\X$ is contravariantly finite and covariantly finite in $\C$;
\item[{\rm (2)}] $M\in\X$ if and only if $\E^k(\X, M)=0$ for any $k\in\{1, 2, ..., n-1\}$;
\item[{\rm (3)}] $M\in\X$ if and only if $\E^k(M,\X)=0$ for any $k\in\{1, 2, ..., n-1\}$.
\end{enumerate}
\end{defn}

Let $\C$ be an extriangulated category with enough projectives and enough injectives.
Given two classes of objects $\X,\Y \subseteq \C$ and an integer $k\geq 1$, the notation $\E^k(\X,\Y) = 0$ will mean that $\E^k(X,Y) = 0$ for every $X\in \X$ and $Y\in \Y$. In the case where $\X = \{ M \}$ or $\Y = \{ N \}$, we shall write $\E^k(M,\Y) = 0$ and $\E^k(\X,N) = 0$, respectively.
The right $k$-th orthogonal complement of $\X$ is defined by
$$\X^{\perp_k} := \{ N \in \C ~|~\E^k(\X,N) = 0 \}.$$
Dually, we have the $k$-th left orthogonal complements ${}^{\perp_k}\Y$.

It is easy to see that $\X$ is an $n$-cluster tilting subcategory of $\C$ if and only if $\X$ is contravariantly finite and covariantly finite in $\C$, and
$$\X = \bigcap\limits_{k=1}^{n-1}{}^{\perp_k}\X=\bigcap\limits_{k=1}^{n-1}{}\X^{\perp_k}.$$

By \cite[Lemma, 2.14]{ZZ}, we know that if $P$ is a projective object, then $\E^k(P,C)=0$ for any $k\geq1$ and $C\in\C$. If $I$ is an injective object, then $\E^k(C,I)=0$ for any $k\geq1$ and $C\in\C$. Hence if
$\X$ is an $n$-cluster tilting subcategory of $\C$, then $\P\subseteq\X$ and $\I\subseteq\X$.

\begin{rem}\label{rem1}
Let $\C$ be an extriangulated category with enough projectives. If $\X$ is a contravariantly finite subcategory in $\C$, then any object $C\in\C$, take a right $\X$-approximation $g\colon X_0\to C$. Since $\C$ has enough projectives, there exists a deflation
$\alpha\colon P\to C$ where $P\in\P\subseteq\X$. By Corollary 3.16 in \emph{\cite{NP}}, we know that $(g,\alpha)\colon X_0\oplus P\to C$ is also a deflation.
Thus there exists an $\E$-triangle $$\xymatrix{B\ar[r]&X_0\oplus P\ar[r]^{\quad(g,~\alpha)}&C\ar@{-->}[r]&}$$
Since $g$ is a right $\X$-approximation of $C$, we have that $(g,\alpha)$ is a right $\X$-approximation of $C$.
Dually, let $\C$ be an extriangulated category with enough injectives. If $\X$ is a covariantly finite subcategory in $\C$. Then for any object $C\in\C$, there exists an $\E$-triangle:
$$\xymatrix{C\ar[r]^f&X\ar[r]&L\ar@{-->}[r]&}$$
where $f$ is a left $\X$-approximation of $C$.
\end{rem}

\section{Main result}

Let $\X$ be a class of objects in an extriangulated category $\C$.  For a nonnegative integer $m\geq 0$,
an $\X$-resolution of $C$ of length $m$ is a complex
$$X_m \to X_{m-1} \to \cdots \to X_1 \to X_0 \to C$$
where $X_k\in\X$ for any integer $0\leq k\leq m$. The above complex is determined by the following
$\E$-triangles:
$$\xymatrix{K_1\ar[r]&X_0\ar[r]&C\ar@{-->}[r]&}$$
$$\xymatrix{K_2\ar[r]&X_1\ar[r]&K_1\ar@{-->}[r]&}$$
$$\vdots$$
$$\xymatrix{K_{n-1}\ar[r]&X_{n-2}\ar[r]&K_{n-2}\ar@{-->}[r]&}$$
$$\xymatrix{X_n\ar[r]&X_{n-1}\ar[r]&K_{n-1}\ar@{-->}[r]&}$$
The \emph{resolution dimension of $C$ with respect to $\X$} (or the \emph{$\X$-resolution dimension} of $C$), denoted $\resdim_{\X}(C)$, is defined as the smallest nonnegative integer $m \geq 0$ such that $C$ has a $\X$-resolution of length $m$. If such $m$ does not exist, we set
$\resdim_{\X}(C) := \infty$.
Dually, we have the concepts of \emph{$\X$-coresolutions of $C$ of length $m$} and of  \emph{coresolution dimension of $C$ with respect to $\X$}, denoted by $\coresdim_{\X}(C)$.

We define
\begin{align*}
\X^\wedge_m & := \{ C \in \C~|~\resdim_{\X}(C) \leq m \}, \\
\X^\vee_m & := \{ C \in \C~|~\coresdim_{\X}(C) \leq m \}.
\end{align*}
In particular, we have $\X^\wedge_0=\X$ and $\X^\vee_0=\X$.
\medskip

Motivated by the definition of $n$-cotorsion pairs in abelian categories \cite[Definition 2.2]{HMP}. We define $n$-cotorsion pairs in extrianglated categories.

\begin{defn}\label{d1}
Let $\C$ be an extriangulated category with enough projectives and enough injectives, and let $\X$ and $\Y$ be two classes of objects of $\C$. We call that $(\X,\Y)$ is a \textbf{left $\bm n$-cotorsion pair} in $\C$ if the following conditions are satisfied:
\begin{enumerate}
\item[{\rm (1)}] $\X$ is closed under direct summands.

\item[{\rm (2)}] $\E^k(\X,\Y) = 0$~ for any $1\leq k \leq n$.

\item[{\rm (3)}] For any object $C \in \C$, there exists an $\E$-triangle
$$\xymatrix{K\ar[r]& X\ar[r]&C\ar@{-->}[r]&}$$
where $X\in \X$ and $K\in\Y^{\wedge}_{n-1}$.
\end{enumerate}
Dually, we can define a right $n$-cotorsion pair. If $(\X,\Y)$ is both a left and right $n$-cotorsion pair in $\C$, we call $(\X,\Y)$ an \textbf{$\bm n$-cotorsion pair} in $\C$.
\end{defn}

Note that when $n=1$, an $n$-cotorsion pair is just
a cotorsion pair in the sense of Nakaoka-Palu, see \cite[Definition 4.1]{NP}.

\begin{exm}
Let $\C$ be an extriangulated category with enough projectives and enough injectives.
It is clear that  both $(\P,\C)$ and $(\C,\I)$ are $n$-cotorsion pair.
We will give more examples of $n$-cotorsion pair in Section 3.

\end{exm}
\begin{lem}\label{y1}
Let $\C$ be an extriangulated category with enough projectives and enough injectives. For any class $\X$ of objects of $\C$, the following holds:
$$\bigcap\limits^{n}_{k = 1}{^{\perp_k}\X}\subseteq {}^{\perp_1}\X^\wedge_{n-1}.$$
\end{lem}

\proof For any $M\in\bigcap\limits^{n}_{k = 1}{^{\perp_k}\X}$, we have $\E^k(M,\X)=0$ for any $1\leq k\leq n$.

Let $N\in\X^\wedge_{n-1}$. Then there exists an $\E$-triangle:
\begin{equation}\label{t1}
\begin{array}{l}
\xymatrix{K_{n-2}\ar[r]& X_{n-1}\ar[r]&N\ar@{-->}[r]&}
\end{array}
\end{equation}
where $X_{n-1}\in\X$ and $K_{n-2}\in\X^{\wedge}_{n-2}$. Apply the functor $\Hom_{\C}(M,-)$ to
the $\E$-triangle (\ref{t1}), we have the following exact sequence:
$$0=\E(M,X_{n-1})\xrightarrow{~~}\E(M,N)\xrightarrow{~\simeq~}\E^2(M, K_{n-2})\xrightarrow{~~}\E^2(M,X_{n-1})=0.$$
Since $\E^k(M,\X)=0$ for any $1\leq k\leq n$, we have $\E(M,N)\simeq\E^2(M,K_{n-2})$.

Since $K_{n-2}\in\X^\wedge_{n-2}$, there exists an $\E$-triangle:
\begin{equation}\label{t2}
\begin{array}{l}
\xymatrix{K_{n-3}\ar[r]& X_{n-2}\ar[r]&K_{n-2}\ar@{-->}[r]&}
\end{array}
\end{equation}
where $X_{n-2}\in\X$ and $K_{n-3}\in\X^{\wedge}_{n-3}$. Apply the functor $\Hom_{\C}(M,-)$ to
the $\E$-triangle (\ref{t2}), we have the following exact sequence:
$$0=\E^2(M,X_{n-2})\xrightarrow{~~}\E^2(M,K_{n-2})\xrightarrow{~\simeq~}\E^3(M, K_{n-3})\xrightarrow{~~}\E^3(M,X_{n-2})=0$$
Since $\E^k(M,\X)=0$ for any $1\leq k\leq n$, we have $\E^2(M,K_{n-2})\simeq\E^3(M,K_{n-3})$.

Inductively, continuing this process, there exists an $\E$-triangle:
\begin{equation}\label{t3}
\begin{array}{l}
\xymatrix{K_0\ar[r]& X_1\ar[r]&K_1\ar@{-->}[r]&}
\end{array}
\end{equation}
where $X_1\in\X$ and $K_0\in\X^{\wedge}_{0}=\X$. Apply the functor $\Hom_{\C}(M,-)$ to
the $\E$-triangle (\ref{t3}), we have the following exact sequence:
$$0=\E^{n-1}(M,X_1)\xrightarrow{~~}\E^{n-1}(M, K_1)\xrightarrow{~\simeq~}\E^n(M, K_{0})\xrightarrow{~~}\E^n(M,X_{1})=0$$
Since $\E^k(M,\X)=0$ for any $1\leq k\leq n$, we have $\E^{n-1}(M,K_1)\simeq\E^n(M,K_{0})$.

Note that $K_0\in\X^{\wedge}_{0}=\X$ and $\E^n(M,\X)=0$, it follows that $$\E(M,N)\simeq \E^n(M, K_{0})=0.$$
This shows that $M\in{}^{\perp_1}\X^\wedge_{n-1}$ and then $\bigcap\limits^{n}_{k = 1}{^{\perp_k}\X}\subseteq {}^{\perp_1}\X^\wedge_{n-1}$.   \qed

\begin{lem}\label{y2}
Let $\C$ be an extriangulated category with enough projectives and enough injectives, and let $\X$ and $\Y$ be two classes of objects of $\C$. Then the following statements are equivalent:
\begin{itemize}
\item[{\rm (1)}] $(\X,\Y)$ is a left $n$-cotorsion pair in $\C$.

\item[{\rm (2)}] $\X = \bigcap\limits_{k=1}^{n}{}^{\perp_k}\Y$ and for any object $C \in \C$ there exists an $\E$-triangle
$$\xymatrix{K\ar[r]& X\ar[r]&C\ar@{-->}[r]&}$$
where $X\in \X$ and $K\in\Y^{\wedge}_{n-1}$.
\end{itemize}
\end{lem}

\proof Note that the implication (2) $\Rightarrow$ (1) is trivial. We show that (1) implies (2). Assume that $(\X,\Y)$ is a left $n$-cotorsion pair in $\C$. By Lemma \ref{t1}, we have the containments
$$\X\subseteq \bigcap_{k = 1}^{n}\,{}^{\perp_k}\Y \subseteq {}^{\perp_1}\Y^\wedge_{n-1}.$$
Thus we only need to prove the remaining containment ${}^{\perp_1}\Y^\wedge_{n-1}\subseteq\X$. For any object $M\in {}^{\perp_1}\Y^\wedge_{n-1}$, there exists an $\E$-triangle
$$\xymatrix{K\ar[r]& X\ar[r]^{g}&C\ar@{-->}[r]&}$$
where $X\in \X$ and $K\in\Y^{\wedge}_{n-1}$. Since $\E(M,\Y^\wedge_{n-1})=0$, the above $\E$-triangle is split.
Hence $g$ is a split epimorphism and then $C$ is a direct summand $X$. It follows that $C\in\X$ implies ${}^{\perp_1}\Y^\wedge_{n-1}\subseteq\X$.  \qed
\medskip

Now we discuss the connection between $n$-cotorsion pairs and $(n+1)$-cluster tilting subcategories.

\begin{thm}\label{main}
Let $\C$ be an extriangulated category with enough projectives and enough injectives.
Then for any subcategory $\X$\hspace{-1.2mm} of\hspace{1mm} $\C$ and any integer $n\geq 1$, the following statements are equivalent:
\begin{enumerate}
\item[{\rm (1)}] $(\X,\X)$ is an $n$-cotorsion pair in $\C$.

\item[{\rm (2)}] $\X$ is an $(n+1)$-cluster tilting subcategory of $\C$.
\end{enumerate}
\end{thm}

\proof (1) $\Rightarrow$ (2). By Lemma \ref{y2} and its dual, we have
$$\X =\bigcap\limits_{k=1}^{n}{}^{\perp_k}\X~\textrm{and}~\X = \bigcap\limits_{k=1}^{n}{}\X^{\perp_k}.$$

For any object $C\in\C$, there exists an $\E$-triangle
$$\xymatrix{K\ar[r]& X\ar[r]^{g}&C\ar@{-->}[r]&}$$
where $X\in \X$ and $K\in\X^{\wedge}_{n-1}$. Apply the functor $\Hom_{\C}(\X,-)$ to
the above $\E$-triangle, we have the following exact sequence:
$$\Hom_{\C}(\X,X)\xrightarrow{~\Hom_{\C}(\X,~g)~}\Hom_{\C}(\X,C)\xrightarrow{~~}\E(\X, K).$$

Since $\X =\bigcap\limits_{k=1}^{n}{}^{\perp_k}\X$, by Lemma \ref{y1}, we have $\X\subseteq {}^{\perp_1}\X^\wedge_{n-1}$ and then $\E(\X,K)=0$. This shows that $g$ is a right $\X$-approximation of $C$, hence $\X$ is a contravariantly finite of $\C$.

Dually, we can show that $\X$ is a covariantly finite subcategory of $\C$.
\smallskip

(2) $\Rightarrow$ (1). Now assume that $\X$ is an $(n+1)$-cluster tilting subcategory of $\C$. Then we have that $\X$ is closed under direct summands and that $\E^k(\X,\X)=0$ for any integer $1\leq k \leq n$.

By Remark \ref{rem1}, for any object $C\in\C$, there exists  an $\E$-triangle:
$$\xymatrix{K_0\ar[r]& X_0\ar[r]^{g_0}&C\ar@{-->}[r]&}$$
where $g_0$ is a right $\X$-approximation of $C$. Apply the functor $\Hom_{\C}(\X,-)$ to
the above $\E$-triangle, we have the following two exact sequences:
$$\Hom_{\C}(\X,X_0)\xrightarrow{~\Hom_{\C}(\X,~g_0)~}\Hom_{\C}(\X,C)\xrightarrow{~~}\E(\X, K_0)\xrightarrow{~~}\E(\X, X_0)=0;$$
$$0=\E^k(\X,X_0)\xrightarrow{~~}\E^k(\X,C)\xrightarrow{~\simeq~}\E^{k+1}(\X,K_0)\xrightarrow{~~}\E^{k+1}(\X,X_0)=0.$$
Since $g_0$ is a right $\X$-approximation of $C$ and $\E^k(\X,\X)=0$ for any $1\leq k\leq n$, we have that
$$\E(\X, K_0)=0~\textrm{and}~\E^{k+1}(\X,K_0)\simeq\E^k(\X,C)~\textrm{for any}~1\leq k\leq n-1.$$

Inductively, continuing this process, there exist the following some $\E$-triangles:
$$\xymatrix{K_m\ar[r]& X_m\ar[r]^{g_m\;\;}&K_{m-1}\ar@{-->}[r]&}$$
where $g_m$ is a right $\X$-approximation of $K_{m-1}$ and $1\leq m\leq n$.
Apply the functor $\Hom_{\C}(\X,-)$ to the above those $\E$-triangles, we obtain the following relations hold:
\begin{align*}
\E(\X,K_n) &= 0,\\
\E^2(\X,K_n) & \simeq \E(\X,K_{n-1}) = 0,\\
\vdots\qquad &\qquad\qquad \vdots\qquad\qquad\qquad \qquad\quad\vdots \\
\E^n(\X,K_n) & \simeq \E^{n-1}(\X,K_{n-1}) \simeq \cdots \simeq\E(\X, K_1) = 0.
\end{align*}
It follows that $K_n\in \bigcap\limits_{k = 1}^{n}\,\X^{\perp_{k}}=\X$ and then $K_0\in\X^{\wedge}_{n-1}$.
\smallskip

This shows that $(\X,\X)$ is a left $n$-cotorsion pair in $\C$.
Dually, we can show that $(\X,\X)$ is a right $n$-cotorsion pair in $\C$.  \qed

\medskip

Now we give some examples illustrating our main result.

\begin{exm}
We revisit Example 5.16 presented in \emph{\cite{LN}}.
Let $\Lambda$ be the self-injective Nakayama algebra given by the following
quiver
$$
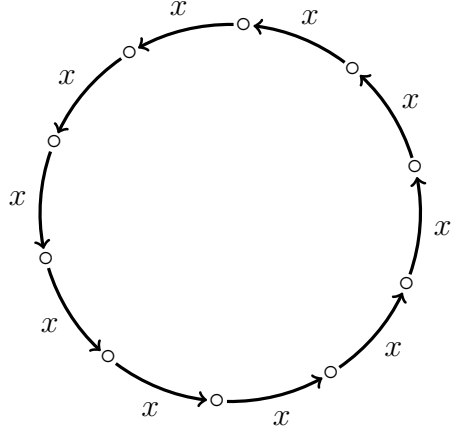
\begin{tikzpicture}
\def \radius {2.5cm}
\def\startDegree{55}
\def\n{10}
\foreach \s in {1,...,\n}
{
  \draw[->, very thick] ({360/\n * (\s - 1)-\startDegree}:\radius)
    arc ({360/\n * (\s - 1)-\startDegree}:{360/\n * (\s - 1)-\startDegree+300/\n}:\radius);
}
\foreach \x in {1,2,...,\n} {
\draw (-36*\x+104:2.8) node {$\\x$} ;
} ;
\foreach \x in {1,2,...,\n} {
\draw (-36*\x+122:\radius) node {$\circ$} ;
} ;
\end{tikzpicture}
$$
with relation $x^4=0$. Then the Auslander-Reiten quiver of the stable category $\underline{{\rm mod}}\Lambda$ of ${\rm mod}\Lambda$ is the following:
$$
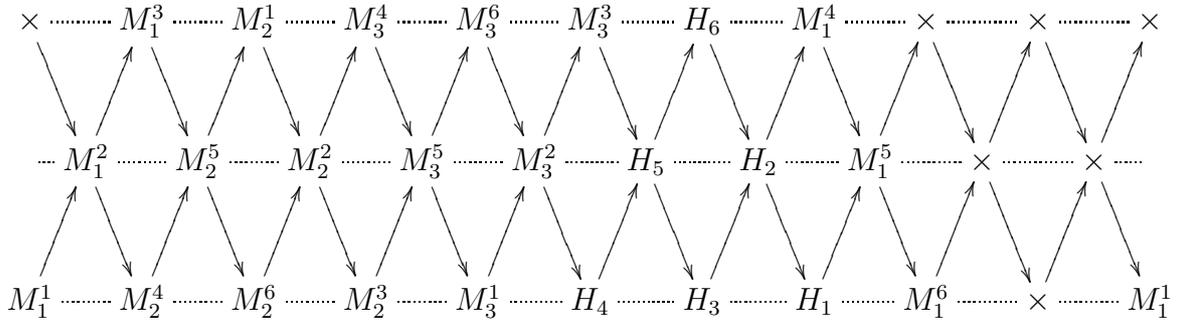
\xymatrix@!@C=-1.2mm@R=10mm{
\times\ar[dr]&&M^3_1\ar@{.}[ll]\ar[dr]&&M^1_2\ar@{.}[ll]\ar[dr]&&M^4_3\ar@{.}[ll]\ar[dr]&&M^6_3\ar@{.}[ll]\ar[dr]&&M^3_3\ar@{.}[ll]\ar[dr]&&H_6\ar[dr]\ar@{.}[ll]&&M_1^4\ar@{.}[ll]\ar[dr]
&&\times\ar@{.}[ll]\ar[dr]&&\times\ar@{.}[ll]\ar[dr]&&\times\ar@{.}[ll]\\
&M_1^2\ar[dr]\ar@{.}[l]\ar[ur]&&M_2^5\ar[dr]\ar[ur]\ar@{.}[ll]&&M^2_2\ar@{.}[ll]\ar[ur]\ar[dr]&&M^5_3\ar@{.}[ll]\ar[ur]\ar[dr]&&M^2_3\ar@{.}[ll]\ar[dr]\ar[ur]&&H_5\ar[ur]\ar[dr]
\ar@{.}[ll]&&H_2\ar@{.}[ll]\ar[dr]\ar[ur]&&M_1^5\ar@{.}[ll]\ar[ur]\ar[dr]
&&\times\ar@{.}[ll]\ar[dr]\ar[ur]&&\times\ar@{.}[ll]\ar[ur]\ar[dr]\ar@{.}[r]&\\
M_1^1\ar[ur]&&M_2^4\ar@{.}[ll]\ar[ur]&&M_2^6\ar@{.}[ll]\ar[ur]&&M^3_2\ar@{.}[ll]\ar[ur]&&M^1_3\ar@{.}[ll]\ar[ur]&&H_4\ar@{.}[ll]\ar[ur]
&&H_3\ar@{.}[ll]\ar[ur]&&H_1\ar@{.}[ll]\ar[ur]&&M_1^6\ar@{.}[ll]\ar[ur]
&&\times\ar@{.}[ll]\ar[ur]&&M_1^1\ar@{.}[ll]}
$$
where the leftmost and rightmost column are identified.

Let $\C$ be the subcategory of the triangulated category $\underline{{\rm mod}}\Lambda$ in which the indecomposable objects are marked by capitals letters.
Since $\C$ is an extension closed subcategory of $\underline{{\rm mod}}\Lambda$, by \emph{\cite[Remark 2.13]{NP}}, we know that $\C$ is an extriangulated category.

Note that $\P={\rm add}( M_1^1\oplus M_1^2\oplus M_1^3)$ (respectively $\I={\rm add}(M_1^4\oplus M_1^5\oplus M_1^6)$) is the subcategory of the projective (respectively, injective) objects. Thus $\C$ has non-trivial projectives
and injectives, which means that it is not triangulated.  It is not exact either, since there
is an inflation $M_1^3\to M_2^5$ which is not monomorphic. In addition, $\C$ has enough projectives.
Indeed, objects $M_2^4,M_2^5,M_2^6$ have deflations from projectives
$$M_1^2\to M_2^4,~M_1^3\to M_2^5,~M_1^3\to M_2^6$$
respectively, and any indecomposable object $N$ outside from $\P\oplus{\rm add}(M_2^4\oplus M_2^5\oplus M_2^6)$
has a deflation $0\to N$. We can also show that $\C$ has enough injectives in a dual manner.

It is straightforward to verify that $\X:={\rm add}(M_1^1\oplus M_1^2\oplus M_1^3\oplus M_1^4\oplus M_1^5\oplus M_1^6)$
is a $4$-cluster tilting subcategory of~ $\C$. By Theorem \ref{main}, we have that $(\X,\X)$ is a $3$-cotorsion pair in $\C$.
\end{exm}

\begin{exm}We revisit Example 3.20 presented in \emph{\cite{V}}. We denote by ``$\circ$" in the Auslander-Reiten quiver the indecomposable objects belong to a subcategory.
Let $\Lambda$ be the algebra given by the following quiver with relations\emph{:}
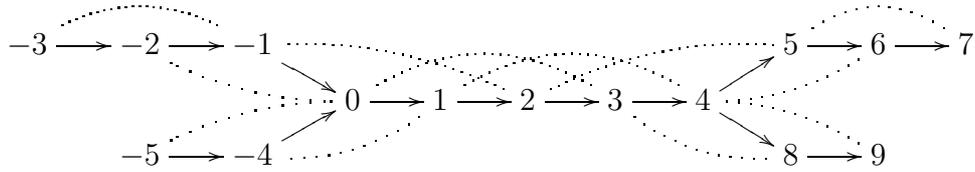
$$\xymatrix@C=0.7cm@R0.2cm{
&\\
-3 \ar[r] \ar@{.}@/^15pt/[rr] &-2 \ar@{.}@/_6pt/[drr]\ar[r] &-1 \ar@{.}@/^8pt/[drrr]\ar[dr] &&&&&&5 \ar[r] \ar@{.}@/^15pt/[rr] &6 \ar[r] &7\\
&&&0 \ar[r] \ar@{.}@/^18pt/[rrr] & 1 \ar[r] \ar@{.}@/^18pt/[rrr] &2 \ar[r] \ar@{.}@/^8pt/[urrr] & 3 \ar[r] \ar@{.}@/_8pt/[drr] &4 \ar[ur] \ar[dr] \ar@{.}@/_6pt/[urr] \ar@{.}@/^6pt/[drr] \\
&-5 \ar[r] \ar@{.}@/^6pt/[urr] &-4 \ar[ur] \ar@{.}@/_8pt/[urr] &&&&&&8 \ar[r] &9
}
$$
There exists a $3$-cluster tilting subcategory $\X$ of~~$\C={\rm mod}\Lambda$\emph{:}
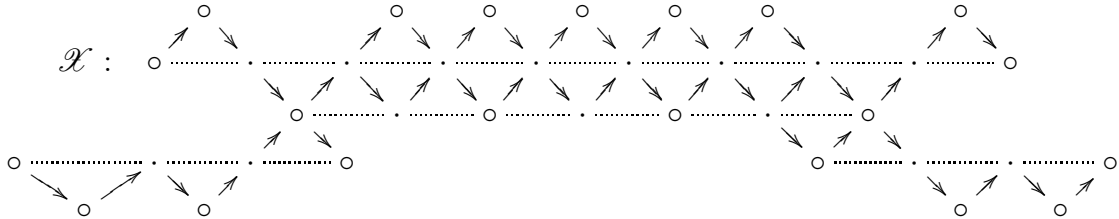
$$\xymatrix@C=0.2cm@R0.2cm{
&&&\circ \ar[dr] &&&&\circ \ar[dr] &&\circ \ar[dr] &&\circ \ar[dr] &&\circ \ar[dr] &&\circ \ar[dr] &&&&\circ \ar[dr]\\
&\X: &\circ \ar@{.}[rr] \ar[ur] &&\cdot \ar[dr] \ar@{.}[rr] &&\cdot \ar[dr] \ar[ur] \ar@{.}[rr] &&\cdot \ar[dr] \ar[ur] \ar@{.}[rr] &&\cdot \ar[dr] \ar[ur] \ar@{.}[rr] &&\cdot \ar[dr] \ar[ur] \ar@{.}[rr] &&\cdot \ar[dr] \ar[ur] \ar@{.}[rr] &&\cdot \ar[dr] \ar@{.}[rr] &&\cdot  \ar[ur] \ar@{.}[rr] &&\circ\\
&&&&&\circ \ar[dr] \ar[ur] \ar@{.}[rr] &&\cdot \ar[ur] \ar@{.}[rr] &&\circ \ar[ur] \ar@{.}[rr] &&\cdot \ar[ur] \ar@{.}[rr] &&\circ \ar[ur] \ar@{.}[rr] &&\cdot \ar[dr] \ar[ur] \ar@{.}[rr] &&\circ \ar[dr] \ar[ur]\\
\circ \ar@{.}[rr] \ar[dr] &&\cdot \ar@{.}[rr] \ar[dr] &&\cdot \ar[ur] \ar@{.}[rr] &&\circ &&&&&&&&&&\circ \ar[ur] \ar@{.}[rr] &&\cdot \ar[dr] \ar@{.}[rr] &&\cdot \ar[dr] \ar@{.}[rr] &&\circ\\
&\circ \ar[ur] &&\circ \ar[ur] &&&&&&&&&&&&&&&&\circ \ar[ur] && \circ \ar[ur]
}
$$
By Theorem \ref{main}, we have that $(\X,\X)$ is a $2$-cotorsion pair in $\C$.
\end{exm}

\begin{exm}
Let $\Lambda$ be a finite-dimensional algebra of global dimension at most $n$.
We denote the Serre functor of $D^b({\rm mod}\Lambda)$
by $\mathbb{S}$, where $D^b({\rm mod}\Lambda)$ is the bounded derived category of ${\rm mod}\Lambda$.
If $\Lambda$ is $n$-representation finite, that is to say, the
module category ${\rm mod}\Lambda$ has an $n$-cluster tilting object, by \emph{\cite[Theorem 1.23]{I}},
we obtain that the subcategory $$\X:={\rm add}\{\mathbb{S}^k\Lambda[-nk]~|~k\in\mathbb{Z}\}$$
of $D^b({\rm mod}\Lambda)$ is $n$-cluster tilting. By Theorem \ref{main},
 we have that $(\X,\X)$ is an $(n-1)$-cotorsion pair in $\C$.
\end{exm}


\begin{thebibliography}{99}
\bibitem{CZZ} W. Chang, P. Zhou, B. Zhu. Cluster subalgebras and cotorsion pairs in Frobenius extriangulated
categories. Algebr. Represent. Theory, 2019, https://doi .org /10 .1007 /s10468 -018 -9811 -7, in press.

\bibitem{HMP} M. Huerta, O. Mendoza, M. A. P\'{e}rez. $n$-Cotorsion pairs. arXiv: 1902.10863, 2019.

\bibitem{HZZ} J. Hu, D. Zhang, P. Zhou. Proper classes and Gorensteinness in extriangulated categories. arXiv:1906.10989, 2019.



\bibitem{I} O. Iyama. Cluster tilting for higher Auslander algebras. Adv. Math. 226: 1-61, 2008.

\bibitem{LN}
Y. Liu, H. Nakaoka.
\newblock Hearts of twin cotorsion pairs on extriangulated categories.
\newblock J. Algebra, 528: 96-149, 2019.


\bibitem{NP}
H. Nakaoka, Y. Palu.
\newblock Extriangulated categories, Hovey twin cotorsion pairs and model structures.
\newblock Cah. Topol. G\'{e}om. Diff\'{e}r. Cat\'{e}g.  60(2): 117-193, 2019.

\bibitem{V}
L. Vaso.
\newblock Gluing of $n$-cluster tilting subcategories for representation-directed algebras.
\newblock arXiv:1805.12180, 2018.


\bibitem{ZZ1}
P. Zhou, B. Zhu.
\newblock Triangulated quotient categories revisited.
\newblock J. Algebra 502: 196-232, 2018.

\bibitem{ZZ} B. Zhu, X. Zhuang. Tilting subcategories in extriangulated categories. arXiv: 1907.00747, 2019.



\end{thebibliography}
\end{document}